\documentclass[11pt]{amsart}
\usepackage{amssymb}
\usepackage{amsmath}
\usepackage{mathtools, float}
\usepackage{amsfonts}
\usepackage{stmaryrd}
\usepackage[english]{babel}
\usepackage[utf8]{inputenc}
\usepackage{enumerate}
\usepackage{graphicx,pinlabel,xcolor, fullpage}
\usepackage{adjustbox}
\usepackage{array,multirow}
\usepackage{float}

\theoremstyle{plain}

\newtheorem{theorem*}{Theorem}

\theoremstyle{definition}

\newtheoremstyle{case}{}{}{}{}{}{:}{ }{}
\theoremstyle{case}

\newcommand{\QED}{\hfill $\blacksquare$}

\begin{document}
\title[Lambda Numbers of Finite Groups]{Lambda Numbers of Finite $p$-Groups}
\author{Mayank Mishra}
\address{Mayank Mishra: \newline
\hspace*{.1in} Plot No.19, Avantipuram, Kalyanpur, Kanpur Nagar, Uttar Pradesh 208 017, India.}
\email{mishra.mayank1903@gmail.com}
\author{Siddhartha Sarkar}
\address{Siddhartha Sarkar: \newline 
\hspace*{.1in} Department of Mathematics, \newline
\hspace*{.1in} Indian Institute of Science Education and Research Bhopal, \newline
\hspace*{.1in} Bhopal Bypass Road, Bhauri, Bhopal 462 066, Madhya Pradesh, India.}
\email{sidhu@iiserb.ac.in}
\urladdr{https://home.iiserb.ac.in/~sidhu/}



\begin{abstract}
An $L(2,1)$-labelling of a finite graph $\Gamma$ is a function that assigns integer values to the vertices $V(\Gamma)$ of $\Gamma$ (colouring of $V(\Gamma)$ by ${\mathbb{Z}}$) so that the absolute difference of two such values is at least $2$ for adjacent vertices and is at least $1$ for vertices which are precisely distance $2$ apart. The lambda number $\lambda(\Gamma)$ of $\Gamma$ measures the least number of integers needed for such a labelling (colouring). A power graph $\Gamma_G$ of a finite group $G$ is a graph with vertex set as the elements of $G$ and two vertices are joined by an edge if and only if one of them is a positive integer power of the other. It is known that $\lambda(\Gamma_G) \geq |G|$ for any finite group. In this paper we show that if $G$ is a finite group of a prime power order, then $\lambda(\Gamma_G) = |G|$ if and only if $G$ is neither cyclic nor a generalized quaternion $2$-group. This settles a partial classification of finite groups achieving the lower bound of lambda number. 
\end{abstract}

\maketitle

\bigskip


\bigskip


\section{\bf Introduction}\label{introsec}

\bigskip

\noindent All graphs considered in this article are assumed to be finite, undirected and simple. Let $\Gamma$ be a graph with vertex set $V(\Gamma)$ and edge set $E(\Gamma)$. For any $a, b \in V(\Gamma)$, the distance $d(a,b)$ is defined to be the length of the shortest path joining the vertices $a$ and $b$ in $\Gamma$. For any $j, k \in {\mathbb{N}}_0 := {\mathbb{N}} \cup \{ 0 \}$, an $L(j,k)$-labelling is defined to be a function $f : V(\Gamma) \rightarrow {\mathbb{Z}}$ that satisfies:

\smallskip

\noindent (i) $|f(a) - f(b)| \geq j$ whenever $d(a,b) = 1$,

\smallskip

\noindent (ii) $|f(a) - f(b)| \geq k$ whenever $d(a,b) = 2$.   

\smallskip

\noindent Therefore an $L(1,0)$-labelling is simply a colouring of the graph $\Gamma$ by non-negative integers. For an $L(j,k)$-labelling $f$ on $\Gamma$, the span of $f$ is defined by
\[
{\mathrm{span}}(f) := \max_{x \in V(\Gamma)} f(x) - \min_{x \in V(\Gamma)} f(x).
\]  
We also define
\[
\lambda_{j,k}(\Gamma) := \min \Big\{ {\mathrm{span}}(f) ~:~ f {\mathrm{~is~an~}} L(j,k){\mathrm{-labelling~on~}} \Gamma \Big\}.
\] 
The {\bf lambda number} of the graph $\Gamma$ is defined by $\lambda(\Gamma) := \lambda_{2,1}(\Gamma)$.

\bigskip

\noindent The $L(j,k)$-labelling of graphs were motivated by the study of the radio channel assignment problem of transmitters through radio or television channels so that the recipients are needed to be separated by certain constraints due to interference (cf. \cite{hal, rob}). One can see \cite{yeh} for a complete survey related to the problem and the references therein.   

\bigskip

\noindent This article is related to the study of lambda number of the power graph $\Gamma_G$ of a finite group $G$. This graph is defined by taking the elements of $G$ as its vertex set, and two distinct vertices $a, b \in G$ are connected by an edge if and only if either $a = b^n$ or $b = a^m$ for some $m, n \in {\mathbb{N}}$. The notion of (directed) power graphs were introduced in \cite{cgs} while studying semigroups. Later on it was shown in a sequel of articles \cite{cgh, cam} that the power graphs of two finite groups are isomorphic if and only if they contain the same number of elements of any given order. In this article we will be focussing on finite groups, and consequently $\Gamma_G$ is always undirected. 

\bigskip

\noindent In \cite{mfw} it is shown that $\lambda(\Gamma_G) \geq |G|$ for any finite group $G$ and the following question is asked:

\bigskip

\noindent {\bf Question.} Classify all finite groups $G$ which satisfy $\lambda(\Gamma_G) = |G|$.      

\bigskip

\noindent Let $p \in {\mathbb{N}}$ be a prime. A finite group $G$ is called a $p$-group if $|G| = p^n$ for some $n \in {\mathbb{N}}$. For such a group $G$, its exponent is the maximum of the orders of elements of $G$. The complete classification of these groups are not known till date. The main problem is the set of parameters those describe the presentations (or even from the viewpoint of representations given over a number field) gets too complicated very quickly as $n = \log_p |G|$ become large. 

\bigskip

\noindent The main theorem of this paper is:

\bigskip

\noindent {\bf Main Theorem.} Let $G$ be a finite $p$-group. Then $\lambda(\Gamma_G) = |G|$ if and only if $G$ is neither cyclic, nor a generalized quaternion $2$-group. 

\bigskip

\noindent Consequently, for a finite $p$-group $G$ of exponent $p^e$ we have
\[
\lambda(\Gamma_G) = 
\begin{cases}
2(p^e - 1) &\mbox{if } G {\mathrm{~is~cyclic}} \\
|G| + 1 &\mbox{if } p=2,~ G {\mathrm{~is~generalized~quaternion}} \\
|G| &\mbox{otherwise }
\end{cases}
\]

\bigskip    

\noindent We discussed some basic group theoretic preliminaries in Section \ref{prelim}. The proof of the main theorem is given in section \ref{proof-main-theorem}. 

\bigskip

\section{\bf Preliminaries}\label{prelim}

\bigskip

\subsection{Definition} Let $G$ be a finite group and $n$ be a positive integer so that $G$ contains an element of order $n$. We define the set $\Lambda(n; G)$ (or simply by $\Lambda(n)$) as   
\[
\Lambda(n;G) := \Big\{ g \in G ~:~ |g| = n \Big\}.
\]
where $|g|$ denote the order of $g$.

\smallskip

\noindent Given any two elements $g_1, g_2 \in G$, we define $g_1 \sim g_2$ if the cyclic subgroups of $G$ generated by $g_1$ and $g_2$ coincide. This is an equivalence relation on $G$. An equivalence class containing $g \in G$ is denoted by $[g]$, and we call it a {\bf cyclic class} of $G$. Then we have 
\[
\Lambda(n;G) = \bigcup_{|g| = n} [g] 
\]
where the union runs over a set of class representatives $g \in G$ with $|g| = n$. In particular, if $n$ is a power of a prime, then the elements of the class $[g]$ induce a complete subgraph of the power graph $\Gamma_G$.

\smallskip

\noindent The class $\{ 1_G \}$ is called the {\bf trivial cyclic class}, where $1_G$ denote the identity element of the group $G$. For any non-trivial cyclic class $[g]$ with $|g| = n$, the class $[g]$ contain $\phi(n)$ number of elements, where $\phi$ denotes the Euler's phi function. For any $n \in {\mathbb{N}}$, we define the $n${\bf -th cyclic class number} ${\mathfrak{m}}(n)$ as the number of cyclic classes that contain the elements of order $n$, i.e.,
\[
{\mathfrak{m}}(n) := {\frac { |\Lambda(n;G)| }{ \phi(n) }}.
\] 

\bigskip

\subsection{Note}\label{triv-adjacency} The $n$-th cyclic class number ${\mathfrak{m}}(n)$ of a group $G$ is equal to the number of cyclic subgroups of $G$ of order $n$. For example, if $G$ is a cyclic group, then ${\mathfrak{m}}(n) = 1$ for every $n$ divides $|G|$.

\bigskip

\noindent Let $G$ be a group and $H, K \leq G$ be subgroups. The subgroup $[H, K] \leq G$ is defined as the subgroup generated by all elements of the form $[h,k] := h^{-1} k^{-1} hk ~(h \in H, k \in K)$. The lower central series of subgroups of $G$ is the descending sequence 
\[
G \geq \gamma_2(G) \geq \gamma_3(G) \geq \dotsc \geq \gamma_i(G) \geq \gamma_{i+1}(G) \geq \dotsc
\]  
of normal subgroups of $G$ given by $\gamma_2(G) := [G, G]$ and $\gamma_{i+1}(G) := [\gamma_i(G), G]$ for every $i \geq 2$. If $G$ is a finite $p$-group, then this descending sequence contains only finitely many non-trivial terms (in general such groups are called {\bf nilpotent}).

\bigskip

\noindent A finite $p$-group of order $p^n$ is said to be of {\bf maximal class} if $\gamma_{n-1}(G) \neq \{ 1_G \}$ and $\gamma_n(G) = \{ 1_G \}$. In this case $G/{\gamma_2(G)} \cong C_p \times C_p$ and $\gamma_i(G)/{\gamma_{i+1}(G)} \cong C_p$ for all $2 \leq i \leq n-1$.

\bigskip

\noindent The following result says a lot more about the class numbers of a finite $p$-group $G$. For a complete proof see \cite[Theorems 1.10, 1.17]{ber2}. 

\bigskip

\subsection{Theorem}\label{class-size} \cite{mil, kula, isa, ber} Let $G$ be a finite $p$-group of exponent $p^e$. Assume that $G$ is not cyclic for an odd prime $p$, and  for $p=2$, it is neither cyclic nor of maximal class. Then:

\smallskip

\noindent (i) ${\mathfrak{m}}(p) \equiv 1 + p~ ({\mathrm{mod}}~ p^2)$.

\smallskip

\noindent (ii) $p \mid {\mathfrak{m}}(p^i)$ for every $2 \leq i \leq e$. 

\bigskip

\subsection{Corollary}\label{thin-class-numbers} Let $G$ be a finite $p$-group of exponent $p^e$. Then ${\mathfrak{m}}(p^i) = 1$ for some $1 \leq i \leq e$ if and only if one of the following occurs:

\smallskip

\noindent (i) $G \cong C_{p^e}$ and ${\mathfrak{m}}(p^j) = 1$ for all $1 \leq j \leq e$, or

\smallskip

\noindent (ii) $p=2$ and $G$ is isomorphic to one of the following $2$-groups:

\smallskip

\noindent (a) dihedral $2$-group 
\[
{\mathbb{D}}_{2^{e+1}} = \Big\langle x, y ~:~ x^{2^e} = 1, y^2 = 1, y^{-1} xy = x^{-1} \Big\rangle, \hspace*{.3in} (e \geq 2)
\] 

\noindent where ${\mathfrak{m}}(2) = 1 + 2^e$ and ${\mathfrak{m}}(2^j) = 1 ~(2 \leq j \leq e)$. 

\smallskip

\noindent (b) generalized quaternion $2$-group
\[
{\mathbb{Q}}_{2^{e+1}} = \Big\langle x, y ~:~ x^{2^e} = 1, x^{2^{e-1}} = y^2, y^{-1} xy = x^{-1} \Big\rangle, \hspace*{.3in} (e \geq 2)
\]

\noindent where ${\mathfrak{m}}(4) = 1 + 2^{e-1}$ and ${\mathfrak{m}}(2^j) = 1$ for all $2 \leq j \leq e$ and $j \neq 2$. 

\smallskip

\noindent (c) semi-dihedral $2$-group
\[
{\mathbb{SD}}_{2^{e+1}} = \Big\langle x, y ~:~ x^{2^e} = 1, y^2 = 1, y^{-1} xy = x^{-1 + 2^{e-1}} \Big\rangle, \hspace*{.3in} (e \geq 3)
\]

\noindent where ${\mathfrak{m}}(2) = 1 + 2^{e-1},~ {\mathfrak{m}}(4) = 1 + 2^{e-2}$ and ${\mathfrak{m}}(2^j) = 1$ for all $3 \leq j \leq e$. 

\bigskip

\noindent {\bf Proof.} From Theorem \ref{class-size}, it follows that if the hypothesis holds then $G$ is either cyclic, or else $G$ is a $2$-group of maximal class. From \cite[Theorem 5.1]{falc}, $G$ must be isomorphic to one of the groups given in (ii)(a)--(c).

\smallskip

\noindent To establish that the condition is necessary and sufficient, we only need to verify the values of ${\mathfrak{m}}(2^j) ~(1 \leq j \leq e)$ as given. Since a cyclic group has a unique subgroup of order $2^j$ for every $1 \leq j \leq e$, the statement (i) follows. Notice that each group given in (ii) has a cyclic maximal subgroup $\langle x \rangle$ of order $2^e$. 

\smallskip

\noindent (ii)(a) For $G = {\mathbb{D}}_{2^{e+1}}$, every element outside the maximal cyclic subgroup $\langle x \rangle \leq G$ has order $2$. Counting the element $x^{2^{e-1}} \in \langle x \rangle$ of order $2$, we have ${\mathfrak{m}}(2) = 1 + 2^e$. For $2 \leq j \leq e$, an element of order $2^j$ is a generator of the cyclic subgroup $\langle x^{2^{e-j}} \rangle$ and hence ${\mathfrak{m}}(2^i) = 1$.

\smallskip

\noindent (ii)(b) The group $G = {\mathbb{Q}}_{2^{e+1}}$ has a unique element of order $2$, namely $x^{2^{e-1}}$. Since every element outside $\langle x \rangle$ has order $4$ and $\phi(4) = 2$, with a similar argument as above, the statement follows. 

\smallskip

\noindent (ii)(c) For $G = {\mathbb{SD}}_{2^{e+1}}$, every element outside $\langle x \rangle$ can be uniquely written as $x^k y ~(0 \leq k \leq 2^e - 1)$. From the given relations it is straightforward to verify
\[
|x^k y| = 
\begin{cases}
2 &\mbox{if } k {\mathrm{~is~even}} \\
4 &\mbox{if } k {\mathrm{~is~odd}} 
\end{cases}.
\]
Counting $x^{2^{e-1}}$, then we have exactly $1 + 2^{e-1}$ of elements of order $2$. Similarly, counting the two elements of order $4$ from $\langle x^{2^{e-2}} \rangle$, we have $2 + 2^{e-1}$ elements of order $4$ and these form a disjoint union of classes of size $\phi(4) = 2$. Now all elements of order $\geq 8$ must belong to $\langle x \rangle$. \QED 

\bigskip  

\noindent We now look at a result from \cite{mfw} and obtain some insights that are going to make the proofs of the next section shorter and easier to understand. 

\bigskip

\subsection{Definition} Let $\Gamma$ be a graph. The {\bf complement graph} $\Gamma^c$ is defined as the graph with the same vertex set $V(\Gamma)$ and two vertices $a, b \in \Gamma^c$ are joined by an edge if and only if they are not adjacent in $\Gamma$.  

\bigskip

\subsection{Theorem}\label{inequality-lambda-number} (\cite[Theorem 3.1]{mfw}) Let $G$ be finite group. Then $\lambda(\Gamma_G) \geq |G|$, with equality if and only if the complement graph $(\Gamma_G \setminus \{ 1_G \})^c$ contains a Hamiltonian path.

\bigskip

\noindent Let $G$ be a finite group $G$ with $\lambda(\Gamma_G) = |G|$. We set $N := |G| - 2$. If $(x_0, x_1, \dotsc, x_N)$ is a Hamiltonian path in $(\Gamma_G \setminus \{ 1_G \})^c$, then the function $f : G \rightarrow {\mathbb{Z}}$ given by $f(1_G) := -2$ and $f(x_i) := i ~(0 \leq i \leq N)$ defines an $L(2,1)$-labelling of $\Gamma_G$ that achieves the lower bound ${\mathrm{span}}(f) = |G|$. 

\bigskip

\noindent On the other hand, if $f : G \rightarrow {\mathbb{Z}}$ is an $L(2,1)$-labelling of $G$ with ${\mathrm{span}}(f) = |G|$, then the same holds for a translation map $f + \alpha$ for any $\alpha \in {\mathbb{Z}}$. Now by definition $f$ is an injective map and from Pigeohole principle it follows that the image ${\mathrm{Im}}(f)$ of the function $f$ contains precisely $|G|+1$ integers with $f(1_G)$ is either the maximum or minimum of ${\mathrm{Im}}(f)$. In case $f(1_G)$ is a maximum of ${\mathrm{Im}}(f)$ we can define an $L(2,1)$-labelling $f_1$ by
\[
f_1(x) = 
\begin{cases}
\min {\mathrm{Im}}(f) - 2 &\mbox{if } x = 1_G \\
f(x) &\mbox{if } x \neq 1_G.
\end{cases}
\]    
Further by applying a translation we can assume that $f_1(1_G) = -2$ and ${\mathrm{Im}}(f) = \{ -2, 0, 1, \dotsc, N \}$. Then the ordered sequence $(f_1^{-1}(0), f_1^{-1}(1), \dotsc, f_1^{-1}(N))$ of vertices in $G \setminus \{ 1_G \}$ is a Hamiltonian path in $(\Gamma_G \setminus \{ 1_G \})^c$. 

\bigskip

\section{\bf Proof of the main theorem}\label{proof-main-theorem}

\bigskip

\noindent We need some preparation before we discuss the proof of the main result. 

\bigskip

\subsection{Definition} Let $G$ be a non-trivial finite group. Two distinct cyclic classes ${\mathcal{C}}_1$ and ${\mathcal{C}}_2$ of $G$ are called {\bf adjacent} if there are elements $g_1 \in {\mathcal{C}}_1$ and $g_2 \in {\mathcal{C}}_2$ so that $g_1$ and $g_2$ are adjacent in the power graph $\Gamma_G$.

\smallskip 

\noindent Clearly, if the cyclic classes ${\mathcal{C}}_1$ and ${\mathcal{C}}_2$ are adjacent, then any pair of elements $x_1 \in {\mathcal{C}}_1$ and $x_2 \in {\mathcal{C}}_2$ are adjacent in $\Gamma_G$. Moreover if the cyclic classes ${\mathcal{C}}_1$ and ${\mathcal{C}}_2$ are adjacent and are represented by elements of the same order, then ${\mathcal{C}}_1 = {\mathcal{C}}_2$.

\bigskip

\subsection{Lemma}\label{lower-hook} ({\bf Lower hook lemma}) Let $G$ be a finite $p$-group and ${\mathcal{U}}, {\mathcal{V}}_1$ and ${\mathcal{V}}_2$ be three classes represented by elements $u, v_1$ and $v_2$ respectively so that $\max \{ |v_1|, |v_2| \} \leq |u|$. If ${\mathcal{U}}$ is adjacent to both of ${\mathcal{V}}_1$ and ${\mathcal{V}}_2$, then ${\mathcal{V}}_1$ is adjacent to ${\mathcal{V}}_2$. In particular, if $|v_1| = |v_2|$, then ${\mathcal{V}}_1 = {\mathcal{V}}_2$.

\begin{figure}[h!]
\centering
\includegraphics[width=0.4\textwidth]{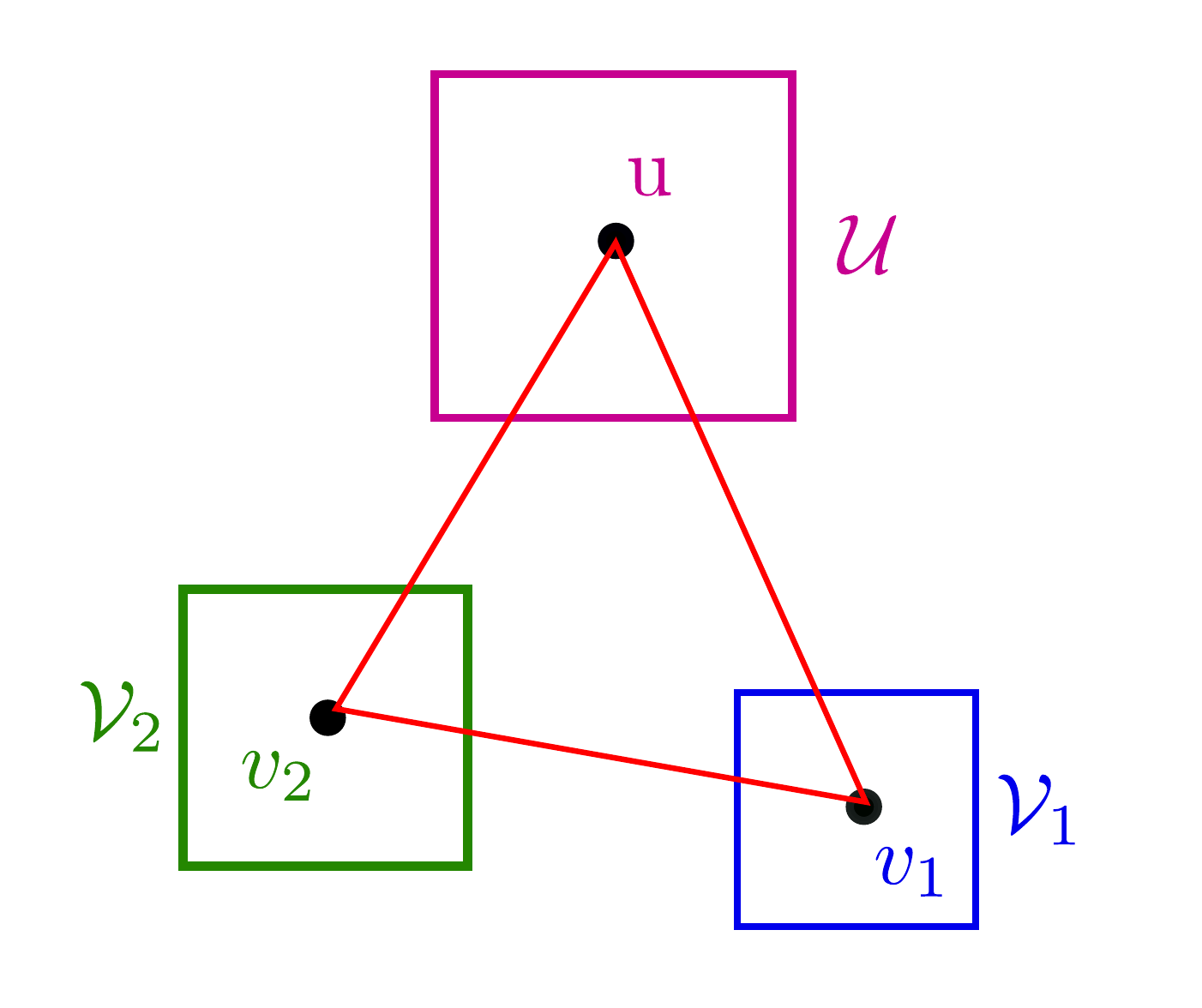}
\caption{}
\label{fig:}
\end{figure}

\bigskip

\noindent {\bf Proof.} Since $u$ is adjacent to $v_1$ and $|v_1| \leq |u|$, we have $\langle v_1 \rangle \subseteq \langle u \rangle$. Similarly, $\langle v_2 \rangle \subseteq \langle u \rangle$. Without loss of generality assume that $|v_1| \leq |v_2|$. Since $\langle u \rangle$ is a cyclic group of prime power order, it follows that $|v_1|$ divides $|v_2|$ and consequently $\langle v_1 \rangle \subseteq \langle v_2 \rangle \subseteq \langle u \rangle$. Now the conclusion follows from the fact that the elements of $\langle u \rangle$ induce a complete subgraph of $\Gamma_G$. \QED

\bigskip

\subsection{Note} The lower hook lemma \ref{lower-hook} is not true for arbitrary finite groups. The cyclic group $C_6$ of order $6$ has four classes ${\mathcal{C}}_1, {\mathcal{C}}_2, {\mathcal{C}}_3$ and ${\mathcal{C}}_6$ represented by elements of orders $1, 2, 3$ and $6$ respectively. Here ${\mathcal{C}}_6$ is adjacent to both of ${\mathcal{C}}_2$ and ${\mathcal{C}}_3$; but the classes ${\mathcal{C}}_2$ and ${\mathcal{C}}_3$ are not adjacent.  

\bigskip

\subsection{Definition}\label{Hamiltonian-complement-graph} Let $\Gamma$ be a graph and ${\mathfrak{M}} := ({\mathcal{M}}_1, {\mathcal{M}}_2, \dotsc, {\mathcal{M}}_r)$ be an ordered collection of subsets of $V(\Gamma)$ of equal size $N$ for some $r \geq 2$. Assume that: 

\smallskip

\noindent (i) each ${\mathcal{M}}_i$ induces a complete subgraph of $\Gamma$, and

\smallskip

\noindent (ii) for any $1 \leq i \neq j \leq r$ and any $a \in {\mathcal{M}}_i,~ b \in {\mathcal{M}}_j$, the vertices $a$ and $b$ are not adjacent in $\Gamma$. 

\smallskip

\noindent We order the elements of ${\mathcal{M}}_j$ as follows: 
\[
{\mathcal{M}}_j = \Big\{ w_{j1}, w_{j2}, \dotsc, w_{jN}  \Big\}, \hspace*{.2in} (j = 1, \dotsc, r).
\]

\begin{figure}[h!]
\centering
\includegraphics[width=0.6\textwidth]{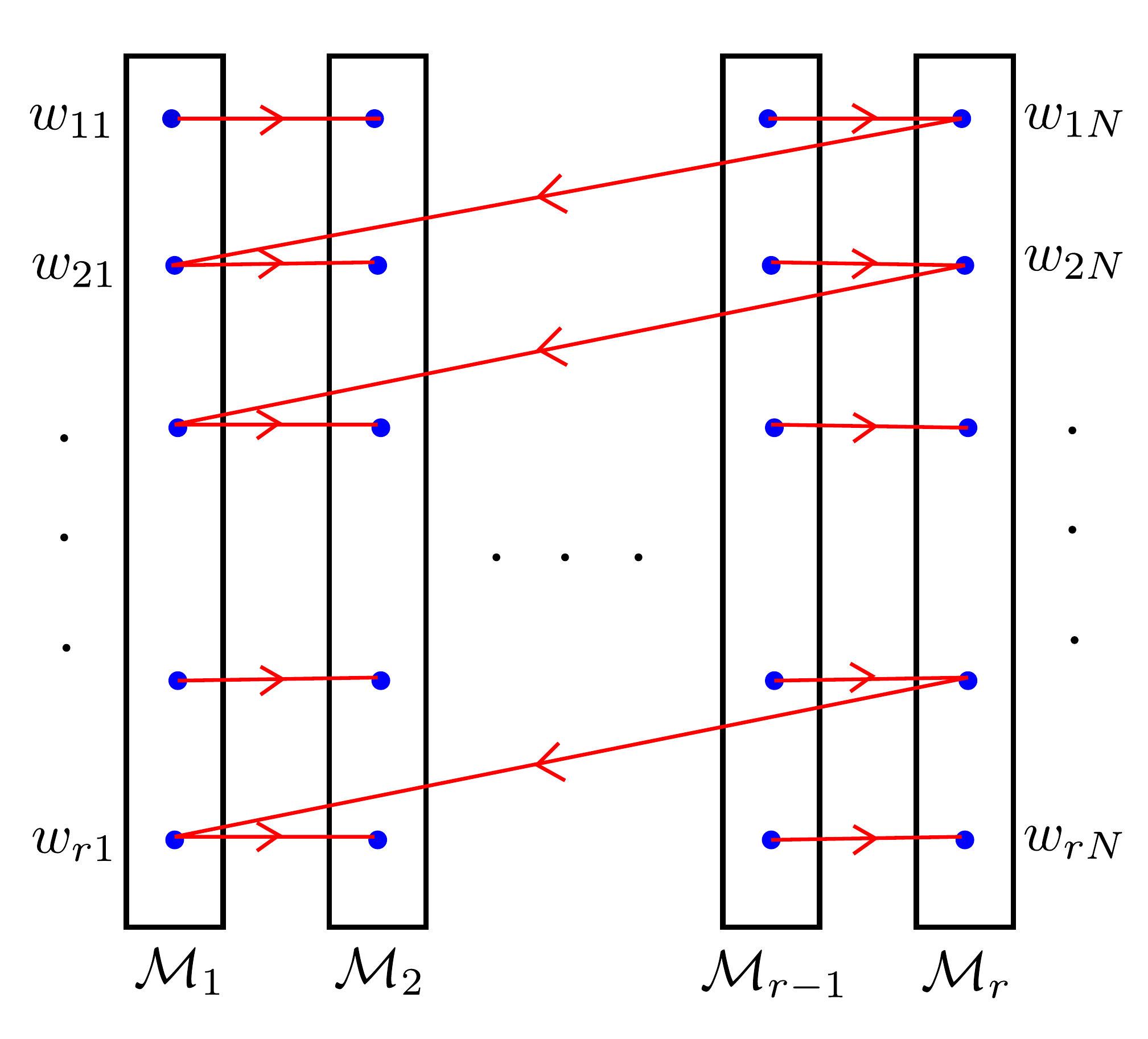}
\caption{${\mathcal{H}}({\mathfrak{M}})$}
\label{fig:Hamiltonian-Plus-fig}
\end{figure}
  
\noindent We define the Hamiltonian path ${\mathcal{H}}({\mathfrak{M}})$ on the set $\cup_{j=1}^{N} {\mathcal{M}}_j$ of vertices (see figure \ref{fig:Hamiltonian-Plus-fig}) as a subgraph of $\Gamma^c$ by:
\[
{\mathcal{H}}({\mathfrak{M}}) := \Big( w_{11}, w_{21}, \dotsc, w_{r-1,1}, w_{r1}, \dotsc, w_{1N}, w_{2N}, \dotsc, w_{r-1,N}, w_{rN} \Big)
\]
We denote by $\partial_0({\mathcal{H}}({\mathfrak{M}})) := w_{11}$ and $\partial_1({\mathcal{H}}({\mathfrak{M}})) := w_{rN}$, the initial and the terminal vertices of the Hamiltonian path ${\mathcal{H}}({\mathfrak{M}})$.

\smallskip

\subsection{Theorem}\label{large-cyclic-class-number} Let $G$ be a finite $p$-group of exponent $p^e$ with ${\mathfrak{m}}(p^i) \geq 2$ for every $1 \leq i \leq e$. Then $\lambda(\Gamma_G) = |G|$. 

\bigskip

\noindent {\bf Proof.} For each $1 \leq i \leq e$, we set $m_i := {\mathfrak{m}}(p^i)$ and let ${\mathcal{C}}^{(i)}_1, {\mathcal{C}}^{(i)}_2, \dotsc, {\mathcal{C}}^{(i)}_{m_i}$ denote the cyclic classes of $G$ containing elements of order $p^i$. Now from lower hook lemma \ref{lower-hook} it follows that, for any $2 \leq i \leq e$ the cyclic class ${\mathcal{C}}^{(i)}_{m_i}$ is non-adjacent to one of the cyclic classes ${\mathcal{C}}^{(i-1)}_1, {\mathcal{C}}^{(i-1)}_2, \dotsc, {\mathcal{C}}^{(i-1)}_{m_{i-1}}$. By relabelling if necessary, we may assume that ${\mathcal{C}}^{(i)}_{m_i}$ is non-adjacent to ${\mathcal{C}}^{(i-1)}_1$.

\bigskip

\noindent Now for each $1 \leq i \leq e$, we set ${\mathfrak{M}}_i := ({\mathcal{C}}^{(i)}_1, {\mathcal{C}}^{(i)}_2, \dotsc, {\mathcal{C}}^{(i)}_{m_i})$ and from the definition \ref{Hamiltonian-complement-graph} we construct the Hamiltonian path ${\mathcal{H}}({\mathfrak{M}}_i)$ in $(\Gamma_G)^c$ on the vertex set $\Lambda(p^i;G)$. Then for every $e \geq i \geq 2$, the vertices $\partial_1({\mathcal{H}}({\mathfrak{M}}_i))$ and $\partial_0({\mathcal{H}}({\mathfrak{M}}_{i-1}))$ are not adjacent in the power graph $\Gamma_G$, and consequently they are connected by an edge $\beta_i$ in the complement graph $(\Gamma_G)^c$. We can now construct the Hamiltonian path ${\mathfrak{H}}$ in $(\Gamma_G \setminus \{ 1_G \})^c$ defined by
\[
{\mathfrak{H}} := \Big( {\mathcal{H}}({\mathfrak{M}}_e), \beta_e, {\mathcal{H}}({\mathfrak{M}}_{e-1}), \beta_{e-1}, {\mathcal{H}}({\mathfrak{M}}_{e-2}), \dotsc, \beta_2, {\mathcal{H}}({\mathfrak{M}}_1) \Big)
\]  
\QED

\bigskip

\noindent {\bf Proof of the main Theorem.} We will first show that for the semi-dihedral $2$-group $G = {\mathbb{SD}}_{2^{e+1}}$ of exponent $2^e ~(e \geq 3)$ we have $\lambda(\Gamma_G) = |G|$ by constructing a Hamiltonian path in $(\Gamma_G \setminus \{ 1_G \})^c$. We are going to use the notations as given in the presentation of $G$ in Lemma \ref{thin-class-numbers}. 

\smallskip

\noindent The cyclic maximal subgroup $\langle x \rangle \leq G$ has order $2^e$ and the elements outside this subgroup are given by $x^k y ~(0 \leq k \leq 2^e -1)$. As noticed earlier, the order of $x^k y$ is $2$ (resp. $4$) if $k$ is even (resp. $k$ is odd). 

\smallskip

\begin{figure}[h!]
\centering
\includegraphics[width=0.8\textwidth]{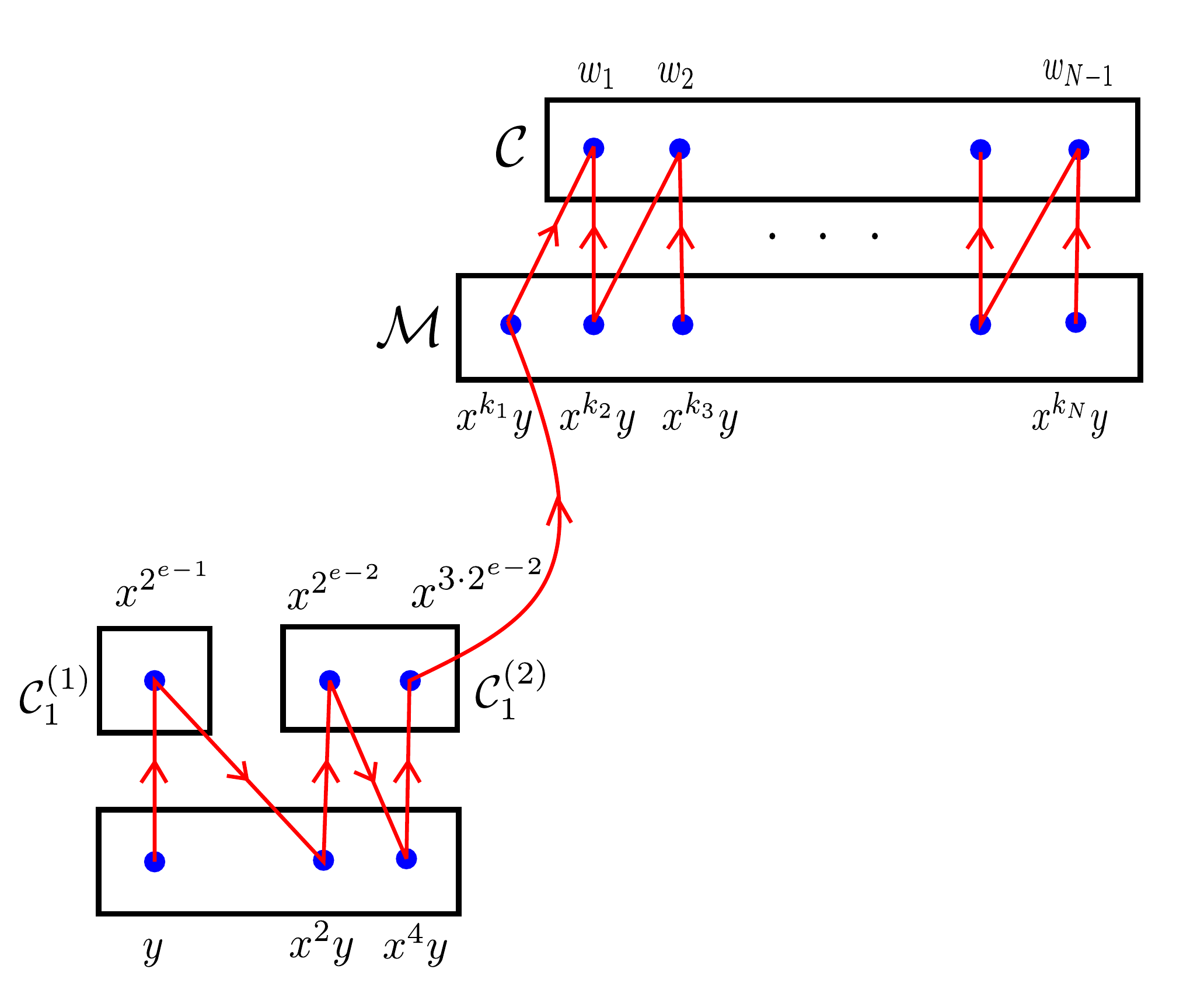}
\caption{}
\label{fig:}
\end{figure}

\noindent We first consider the cyclic classes of elements of order $2$ and $4$ from the cyclic subgroup $\langle x \rangle$; namely ${\mathcal{C}}^{(1)}_1 = \{ x^{2^{e-1}} \} \subseteq \Lambda(2;G)$ and ${\mathcal{C}}^{(2)}_1 = \{ x^{2^{e-2}}, x^{3 \cdot 2^{e-2}} \} \subseteq \Lambda(4;G)$. Now we consider the sequence
\[
T_1 := \Big( y,~ x^{2^{e-1}},~ x^2 y,~ x^{2^{e-2}},~ x^4 y,~ x^{3 \cdot 2^{e-2}} \Big).
\]  
Clearly this forms a path in the complement graph $(\Gamma_G)^c$. Now we fix an arbitrary labelling of the set $(G \setminus \langle x \rangle) \setminus \{ y, x^2 y, x^4 y \}$ as follows 
\[
{\mathcal{M}} = \Big\{ x^{k_i} y ~:~ i = 1, 2, \dotsc, N \Big\} 
\]
where $N = 2^e - 3$. Notice that the number of elements in ${\mathcal{C}} = \langle x \rangle \setminus \{ 1_G, x^{2^{e-1}}, x^{2^{e-2}}, x^{2^{e-1} + 2^{e-2}} \}$ is $2^e - 4 = N - 1$.

\smallskip

\noindent We claim that for any pair of elements $a \in {\mathcal{C}}$ and $b \in {\mathcal{M}}$, the vertices $a$ and $b$ are not adjacent in $\Gamma_G$ (and hence they are adjacent in the complement $(\Gamma_G)^c$). To see this, notice that $a$ is a power of $x$ and has order $\geq 8$, whereas $b$ has order either $2$ or $4$. So if they are adjacent, then $b = a^m$ for some positive integer $m$ and consequently, $b \in \langle x \rangle$, a contradiction.  

\smallskip

\noindent Now setting up an arbitrary enumeration of the elements of ${\mathcal{C}} = \{ w_1, w_2, \dotsc, w_{N-1} \}$ we have the following path in $(\Gamma_G)^c$
\[
T_2 = \Big( x^{k_1} y,~ w_1,~ x^{k_2} y,~ w_2, \dotsc,~ x^{k_{N-1}} y,~ w_{N-1}, x^{k_N} y \Big).
\]

\smallskip

\noindent Now a similar argument as above shows that the elements $x^{k_1} y$ and $x^{2^{e-1} + 2^{e-2}}$ are not adjacent in $\Gamma_G$. Hence there exists an edge $\beta$ in $(\Gamma_G)^c$ joining these two vertices. Consequently, the path $(T_1, \beta, T_2)$ is a Hamiltonian path in $(\Gamma_G \setminus \{ 1_G \})^c$. This proves that $\lambda(\Gamma_G) = |G|$.   

\smallskip

\noindent Now we are ready to prove the final statement. Let $G$ be a finite $p$-group of exponent $p^e$. Assume that $\lambda(\Gamma_G) = |G|$. Since $\lambda(\Gamma_{C_{p^e}}) = 2(p^e - 1)$ (since $\Gamma_{C_{p^e}}$ is a complete graph) and $\lambda({\mathbb{SD}}_{2^{e+1}}) = 2^{e+1} + 1$ (see \cite[Example 3.3]{mfw}), $G$ cannot be these groups.

\smallskip

\noindent Conversely, suppose $G$ is neither cyclic, nor a generalized quaternion $2$-group. If ${\mathfrak{m}}(p^i) \geq 2$ for every $1 \leq i \leq e$, then from Theorem \ref{large-cyclic-class-number}, we are done. Hence using Corollary \ref{thin-class-numbers}, we have $p=2$ and $G$ is isomorphic to one of the groups ${\mathbb{D}}_{2^{e+1}}$ or ${\mathbb{SD}}_{2^{e+1}}$. If $G \cong {\mathbb{D}}_{2^{e+1}}$, then from \cite[Example 3.2]{mfw}, the statement follows. Finally if $G \cong {\mathbb{SD}}_{2^{e+1}}$, then we have shown above that $\lambda(\Gamma_G) = |G|$. \QED  

\vspace*{.3in} 
 
\bibliographystyle{plain} 
\bibliography{LambdaNumberFGV1}


\vspace*{.3in}

\end{document}